\documentclass[12pt, a4paper]{article}

\usepackage[usenames]{color}

\usepackage{amssymb,amsthm,amsmath,tocloft,longtable}
\title{Bounds on Multiplicities of Laplace Operator Eigenvalues on Surfaces}

\author{Aleksandr Berdnikov\thanks{
Independent University of Moscow, Bolshoy Vlasyevskiy per.~11, 119002, Moscow, Russia;
\newline 
Department of Mathematics, Massachusetts Institute of Technology, 77 Massachusetts Ave, Cambridge, MA 02139, USA; 
\newline
National Research University Higher School of Economics,
Vavilova Str.~7, 117312, Moscow, Russia.
\newline
{\tt aberdnik@mit.edu}}}
\date{2015}

\begin{document}

\maketitle

\newtheorem{lema}{Lemma}
\newtheorem{prop}{Proposition}
\newtheorem{teo}{Theorem}

\section{Introduction}

Hoffmann-Ostenhof, Mihor and Nadirashvili proved in the paper~\cite{orig} several bounds on multiplicities of eigenvalues of the Laplacian operator on a simply-connected planar domain with Dirichlet boundary condition, which were proved for compact surfaces by Nadirashvili in~\cite{porig}. In the present paper we generalize these bounds. It was also claimed in~\cite{orig} that the original proof works in the case of non-simply-connected domains after taking minor remarks into account. We show that such a modified proof contains a gap.

We use the following notation

\medskip

\begin{tabular}{|p{2.3 cm}|p{12 cm}|}
\hline
$\Delta$ & Laplace operator \\
\hline
$\lambda_k$ & $k$-th Laplacian eigenvalue, $0\leq \lambda_1\leq \dots \leq \lambda_k \leq \dots$\\
\hline
$\sigma_k$ & $k$-th Steklov eigenvalue, $0\leq \sigma_1\leq \dots \leq \sigma_k \leq \dots$\\
\hline
$U_k$ & Laplacian (Steklov) eigenspace corresponding  to $\lambda_k$ ($\sigma_k$)\\
\hline
$m(\lambda_k), m(\sigma_k)$ or $m(k)$ & Multiplicity of $k$-th eigenvalue, $\dim (U_k)$\\
\hline
$\mathcal{N}(f)$ & Nodal graph of a function $f$\\
\hline
$\mu (f)$ & Number of connected components of $M\setminus \mathcal{N}(f)$\\
\hline
$Op(X)$ & A sufficiently small neighborhood of the set $X$\\
\hline
$Int(X)$ & Interior of $X$\\
\hline
\vspace{-3mm}
$(\partial M)^i$ & $i$-th connected component of the boundary $\partial M$ of a surface $M$\\
\hline
\vspace{-3mm}
$\bar{M}$ & \vspace{-3mm} A surface obtained from $M$ by contracting each $(\partial M)^i$ to a point\\
\hline
$pt^i$ & Image of $(\partial M)^i$ in $\bar{M}$\\
\hline
\end{tabular}

\medskip

There is a number of known bounds for the multiplicities of Laplacian and Steklov eigenvalues. The classical estimate 
\begin{equation}
\label{c}
m_k\leq 2k-2\chi(M)-2b_0(\partial M)+3
\end{equation}
where $b_0(\partial M)$ is the number of connected components of $\partial M$, was proven in various generality and in different cases in~\cite{besso},~\cite{polt},~etc. These estimates can be significantly improved in the case of the Steklov problem: in~\cite{polt} it was shown that in the Steklov case the following estimate holds
\begin{equation}
\label{st}
m_k\leq k-2\chi(M)+3
\end{equation}
and under some assumptions the inequality~\eqref{st} is strict, and even more, for each surface $M$ there is a constant $C(M)$ such that $m_k(M)\leq C(M)$. In the case of Laplacian on the closed surfaces $M$ with $\chi(M)<0$ the estimate~\eqref{c} is improved in~\cite{porig} by two: 
\begin{equation}
\label{cc}
m_k\leq 2k-2\chi(M)+1.
\end{equation}
Later this result was generalized in~\cite{orig} to the case of a simply-connected planar domains with the Dirichlet boundary condition. In the present paper we prove that $m_k\leq 2k-2\chi(M)-2b_0(\partial M)+1$ in the case of surfaces $M$ with boundary such that $\bar{\chi}{:}=\chi(M)+b_0(\partial M)<0$ in both Laplacian and Steklov problem. This result generalises~\eqref{cc} which improves in the case of Laplacian the known estimate~\eqref{c} by two, but in the Steklov case it is better than both~\eqref{c} and~\eqref{st} only for $k\leq 2b_0(\partial M)$.

There in a number of better estimates for muliplicities of few first eigenvalues, some of which were proven to be sharp, see for example~\cite{jammes}.

\bigskip

The main result of the paper~\cite{orig} is
\medskip

{\bfseries Theorem A}. {\itshape Let $k\geq 3$. Then the multiplicity of the $k$-th eigenvalue $\lambda_k$ for the Dirichlet problem on a simply-connected planar domain $D$ satisfies}
$$m(\lambda_k) \leq 2k-3.$$

Theorem A immediately follows from the following results.
\medskip

{\bfseries Theorem B}~\cite{orig}. {\itshape Let $U$ be a linear subspace of an eigenspace $U_k$ and let $1<l\in\mathbb{N}$ be an integer such that for each $f\in U$ we have $\mu (f)\leq l$. Then $\dim (U)\leq max(3,2l-3)$}

\medskip

{\bfseries Courant nodal domain theorem} \cite{courant}. {\itshape For each function $f$ in the eigenspace $U_k$ we have $\mu (f)\leq k$}.
\medskip

The main results of the present paper are the following analogs of theorems A and B valid for a surface $M$ with boundary.

\begin{teo}
\label{t1}
Let $\lambda$ be a real number and $U$ be a linear space of functions on a surface $M$ with a boundary $\partial M$ consisting of $b$ connected components. Suppose that every $f\in U$ is a Laplacian eigenfunction $\Delta f=\lambda f$. Let $l=\sup \{ \mu (f)\mid f\in U\}$ and suppose $\chi(M)+b<0$. Then $\dim (U)\leq 2l-2(\chi(M)+b)+1$.
\end{teo}

Remark that no boundary conditions are imposed on $f\in U$. 

According to Courant theorem the is an estimate $\sup \{ \mu (f)\mid f\in U_k\}\leq k$, hence the the theorem~\ref{t1} can be applied to the eigenspace $U_k$ with $l\leq k$. As a result we get the following theorem.

\begin{teo}
\label{t2}
Let $M$ be the surface with $\chi(M)+b<0$. Then the multiplicity $m(k)$ of $k$-th Laplace or Steklov eigenvalue (for any given boundary conditions) satisfies
$$m(k)\leq 2k-2(\chi (M)+b)+1.$$
\end{teo}

\section{Sketch of the Proof}

First of all we introduce "star fibrations" $E_M(n)$ over any surface $M$. This helps to formalize the argumentation used in~\cite{orig},~\cite{porig} and the present paper. This fibrations parametrize certain types of singularities that the Laplacian eigenfunctions may possess.

Assuming that the multiplicity $m(k)$ is bigger than the estimate we are aimed to prove, we find for each point $x \in M$ an eigenfunction $f_x$ possessing the considered type of singularity at the point $x$, and providing therefore an element $s(x)$ of a fibre $F_M$ of $E_M$. We prove that under the assumptions taken, $s(x)$ is uniquely defined for any point $x\in M$ and depends on $x$ smoothly. Thus, we obtain a smooth section $s$ of $E_M$.

If the surface $M$ is closed, then the existence of the section $s$ of $E_M$ implies that the Euler class of the orienting cover of $M$ is zero. Since $\chi (M)<0$ was supposed, we obtain a contradiction and the assumption of high multiplicity $m(k)$ was wrong. That is how the proof of~\cite{porig} is organized.

If, however, $M$ is not closed, then $H^2(M)=0$ and Euler class is zero for any such surface; $E_M$ has actually got a lot of sections and no contradiction arise yet. To deal with this, Nadirashvili at. al. study the behaviour of $s$ near the boundary $\partial M$ and obtain certain homotopical restriction on $s$. If $M=D^2$ then no section satisfying this restriction on $\partial D^2$ can be extended into the interior of the disc $D^2$ and that is the desired contradiction with the existence of $s$. However, if $M$ has more than one hole, there are a lot of sections satisfying that restriction. That prevents finishing the proof (despite claims of~\cite{orig}).

The described problem can be avoided if $\chi (M)+b<0$ (namely, if $M$ possesses handles). In this case we obtain the homotopical restriction on $s$ which implies, that the section $s$ can actually be extended onto $\bar{M}$ (the closed surface obtained from $M$ by glueing its holes with discs). Then we obtain the desired contradiction just like in~\cite{porig} and conclude the proof.

\section{Star Fibration}

\label{star}

It is useful to introduce a notion of a star fibration on a riemannian surface for discussing methods of~\cite{orig} used is the present paper. The motivation is provided by the following theorem.

{\bfseries Bers theorem} \cite{bers}. {\itshape For a Laplacain eigenfunction $f$ and $x_0\in M$ there exists an integer $n\geq 0$ and real numbers $A, B\in\mathbb{R}$ such that in polar coordinates $(r,\theta)$ centered at $x_0$ the following formula holds}
$$
f(x)=r^n(A\sin (n\theta)+B\cos (n\theta))+O(r^{n+1}).
$$
This means that the nodal graph of the eigenfunction $f$ is diffeomorphic near $x_0$ to $2n$ rays in $\mathbb{R}^2$ emitting from $0$ at equal angles between the adjacent lines. We say in this situation that the function $f$ has a zero of order $n$ at $x_0$. 

Bers theorem makes it natural to consider a star fibration $E_M(2n)$ over $Int(M)$ defined as follows. The fiber $F_x(2n)$ at a point $x$ of the fibration $E_M(2n)$ consists of all $2n$-stars in the tangent space $T_xM$, i.e. configurations of $2n$ rays in $T_xM$ with equal angles between adjacent lines. Formally, let $SP^{k}(STM)$ denote a fiberwise symmetric product of $k$ copies of the spherisation $STM$ of the tangent bundle and $\rho_\alpha$ denote the rotation by $\alpha$ in $STM$ in some locally chosen direction. Then
$$E_M(2n)=\{(l,\rho_{2\pi \frac{1}{2n}}(l),...,\rho_{2\pi \frac{2n-1}{2n}}(l))\mid l\in STM\}\subset SP^{2n}(STM).$$
It is clear that the choice of the direction of rotation $\rho_\alpha$ does not affect the result.

In these terms the Bers theorem states that the nodal graph $\mathcal{N}(f)$ defines at each its vertex $x$ of degree $2n$ an element of $F_x(2n)$ which we denote by $s(\mathcal{N}(f),x)$. Besides the true nodal graphs $\mathcal{N}(f)$ we consider graphs $\mathcal{N}'\subset \bar{M}$ which are diffeomorphic to nodal graphs $\mathcal{N}(f)$ of some eigenfunctions $f$. Clearly, it is still possible to consider a star $s(\mathcal{N}',x)$ for any vertex $x\in \mathcal{N}'$. 

We will need the following technical lemma about $s$. 

\begin{prop}
Let $\lambda$ be a real number and let $U$ be a finite-dimensional space of functions $f\in U$ satisfying $\Delta (f)=\lambda f$. For every $x\in Int(M)$ consider the subspace $[f_x]\subset U$ of functions $f_x$ with a zero of order at least $n$ at $x$. Suppose that for any $x\in Int(M)$ the subspace $[f_x]$ is of dimension 1 (i. e. nonzero $f_x$ are defined up to $\mathbb{R}^*$), and that the order of zero $x$ of $f_x$ is precisely $n$. Define a section $\sigma (x)=s(\mathcal{N}(f_x),x)$. Then this section $\sigma \in \Gamma (E_M(2n))$ is smooth.
\end{prop}

{\bfseries Proof}. Consider local coordinates $x=(x^1, x^2)$. Then functions $f_x$ can be defined by a system of linear equations in $U$ depending on $x$. This system requires the patrial deravatives of $f_x$ of order $<n$ to vanish:
\begin{equation}
\label{system}
\left\{
\begin{array}{lllll}
f_x(x)=0,& & & & \\
\frac{\partial}{\partial x^1}f_x(x)=0,& &\frac{\partial}{\partial x^2}f_x(x)=0,& \\
\dots & & \\
(\frac{\partial}{\partial x^1})^n f_x(x)=0,& \dots & (\frac{\partial}{\partial x^1})^i (\frac{\partial}{\partial x^2})^{n-i} f_x(x)=0, & \dots &(\frac{\partial}{\partial x^2})^n f_x(x)=0.\\
\end{array}
\right.
\end{equation}

Note that the system \eqref{system} is smooth and of constant rank since its kernel $[f_x]$ is one-dimensional at any $x$. That means that $[f_x]$ depends smoothly on $x$ and we can take smooth represenatives $f_x$. Let $f_x=\sum c_i(x)f_i$ for a basis $\{f_i\}$ of $U$, then $c_i(x)$ are smooth. For any fixed $x\in Int(M)$ the function $f_x$ is a function on $M$. Consider the $n$-th term $f^{(n)}_x$ of its Taylor series at $x$. Since $f_x=\sum c_i(x)f_i$ and both $f_i$ and $c_i(x)$ depend smoothly on all of their arguments, then the term $f^{(n)}_x$ depends on $x$ smoothly as well. Since it is the first non-zero term at $x$, it defines the section $\sigma (x)$. This implies that $\sigma (x)$ depends smoothly on $x$ as well. $\square$

\section{A Gap in the Paper~\cite{orig}}

We discuss in this section the topological methods from the paper~\cite{orig} and show that the proof from the paper~\cite{orig} contains a gap in the case of non-simply-connecned domains.

\begin{prop}
{\upshape \cite{orig}} If $m(\lambda_k)\geq 2k-2$ then for any $x\in Int(M)$ there exists a unique (up to $\mathbb{R}^*$) eigenfunction $f_x\in U_{k}$ such that nodal graph $\mathcal{N}(f_x)$ has a single vertex $x$ of degree $2k-2$ and each component $(\partial M)^i$ of the boundary intersects at most two edges of $\mathcal{N}(f_x)$. 
\end{prop}

Proposition 2 means that the graph $\mathcal{N}(f_x)$ in $\bar M$ has $x$ as a single vertex, where $\bar{M}$ is the space obtained by contracting each $(\partial M)^i$ into a point $pt^i$. It follows from Propositions 1 and 2 that if $m(\lambda_k)\geq 2k-2$ then the smooth section $\sigma (x)=s(\mathcal{N}(f_x),x)$ of $E(2k-2)$ is defined on $Int(M)$.

The next claim we cite (\cite{orig} p. 1185, (1), (2)) can be rewritten in the language of star fibrations in the following way. Take a trivialization of $E_M(2k-2)$ over a certain neighborhood $Op((\partial M)^i)$ induced by the representation 
$$Op((\partial M)^i)\cong (\partial M)^i\times [0,1)\cong S^1\times [0,1).$$ 
Then the Claim (2) from~\cite{orig} states that the homotopy class of the section $\sigma |_{(\partial M)^i\times\varepsilon}$ differs from the trivial one by rotation in the positive direction in the following sence. 

Consider the coordinates $(x^1, x^2)$ provided by the decomposition $$Op((\partial M)^i)\cong (\mathbb{R}/\mathbb{Z})\times [0,1).$$ Consider a loop $\gamma (t)=(t, \varepsilon), t\in [0, 1]$. Then any section of \mbox{$E_M(2k-2)|_{\gamma (\mathbb{R})}$} is homotopic to the section $\rho (vt) (s_0)$ where $s_0$ is a constant (in chosen coordinates) section and $\rho (vt)$ is a rotation 
$$\left(\begin{array}{cc}
\cos (vt) & -\sin (vt) \\
\sin (vt) & \cos (vt) \\
\end{array}
\right)$$
in $T_{\gamma (t)}M$ for some $v\in \mathbb{R}$ such that $(2k-2)v \in \mathbb{Z}$, $v$ is the speed of rotation. The claim (2) states that for $\sigma (\gamma (t))$ the corresponding $v$ is non-negative.

It follows from claim (2) that the assumption $m(\lambda_k)>2k-2$ leads to the contradiction in the simply-connected case $M=D^2$. Indeed, according to the claim (2), the section $\sigma$ near $\partial M\cong S^1$ provides non-negative element $h(s)\in \pi_1(F(2k-2))\cong \mathbb{Z}$ where $E(2n)|_{\gamma}$ is trivialized, as above, by outward normal trivialization of $T(Op(\partial M^i))$. Then when we pass to the standard trivialization of $T(D^2)$ the element $h(s)$ increases by $2k-2$ (it is the class of the outward normal vector field in $\pi_1(F(2k-2))$ with respect to the standard trivialization). But the restriction of the section $\sigma$ on any loop in $D^2$ should vanish in $\pi_1(F(2k-2))$ since all loops in $D^2$ are contractible and the section $\sigma$ is defined on the whole $D^2$.

Let us now explain the gap in the proof from the paper~\cite{orig} in the case of a non-simply-connected planar domains. It is remarked in~\cite{orig} that this case can be treated in the analogous way. It is proposed to make several cuts in order to make a disc from $M$ and then apply the same arguments to the loop going along the new boundary. However, here we have a problem caused by the fact that the inner holes' trivialisation provides an addition of a {\it negative} class in $\pi_1(F(2k-2))$ in contrast to the {\itshape positive} class provided by the outer component of the boundary. Hence there is no more any reason for the resulting obstruction we consider to be non-zero. Hence there is no contradiction and the assumption $m(\lambda_k)>2k-2$ might be true.

\section{Proof of Theorem 1}

\label{plan}

We start by proving a statement analogous to Proposition 2.

\begin{lema}

Let the assumptions of Theorem 1 hold, i. e. let $\lambda \in \mathbb{R}$ be a real number, $U$ be a vector space of smooth functions on $M$ satisfying $\Delta (f)=\lambda f$ and let \mbox{$l=\sup \{\mu(\mathcal{N}(f)|f\in U\}$}. Suppose Theorem 1 does not hold, i. e. \mbox{$\dim (U)\geq 2l-2(\chi (M)+b)+2=2l-2\chi (\bar{M})+2$}. Denote $l-\chi (\bar{M})+1$ by $n$. For each $x\in Int(M)$ consider the set $U_n\subset U$, consisting of eigenfunctions $f_x$ whose nodal graph $\mathcal{N}_f$ contains a vertex $x$ of degree at least $2n$. 

Then $\dim (U_n)\geq 1$. Moreover,
\newline
\begin{tabular}{p{14.5 cm}}
$\bullet$ any $f_{x}$ has a nodal graph with a unique vertex in $\bar{M}$; note that if two rays have an endpoint $pt^i$, we do not consider $pt^i$ as a new vertex $x$ in $\bar{M}$; \\
$\bullet$ $\deg_{\mathcal{N}(f_x)}(x)=2l-2\chi (\bar{M})+2$;\\
$\bullet$ $\dim (U_n)\leq 2$; \\
$\bullet$ faces of $\mathcal{N}(f_{x})$ are homeomorphic to $D^2$.\\
\end{tabular} 
\end{lema}

The proof is postponed till the section~\ref{l1}.

\medskip

Lemma 1 implies that under the assumptions of Theorem 1 for each point $x$ there is either unique (up to $\mathbb{R}^*$) eigenfunction $f_x$ whose nodal graph $\mathcal{N}_{f_x}$ has the vertex $x$ of degree $\geq 2l-2\chi (\bar{M})+2$, or a 2-dimentional space $U_n$ of such functions and $P(U_n)\sim S^1$. We prove that the case $\dim (U_n)=2$ described above is impossible.

\begin{lema}
Consider a $S^1$-parameterized continuous family of graphs $\mathcal{N}(t)$ in $\bar{M}$ isotopic to a family of nodal graphs. Suppose that  

\begin{tabular}{p{14.5 cm}}
$\bullet$ $\chi(\bar{M})<0$;\\
$\bullet$ for all $t\in S^1$ the graph $\mathcal{N}(t)$ has the only vertex $x$;\\
$\bullet$ the degree $\deg_{\mathcal{N}(t)}(x)=2n$ is constant;\\
$\bullet$ at least one of edges of $\mathcal{N}(t)$ provides a non-contractible loop $l\in \pi_1(\bar{M},x)$.\\
\end{tabular}

Then induced loop $s(\mathcal{N}(t),x)$ is contractible in the fiber $F_x(2n)$.
\end{lema}

The proof is postponed till the section~\ref{l2}.

\medskip

Lemma 2 excludes the case $\dim (U_n)=2$ due to the following argument. Let us check that Lemma 2 applies to the graphs from Lemma 1, i. e. graphs $\mathcal{N}(f)$ for $f\in U_n$ have non-contractible edges. Indeed, the nodal graph $\mathcal{N}(f)$ provides cellular decomposition of $\bar{M}$, hence the edges of $\mathcal{N}(f)$ generates $H_1(\bar{M})$. Since $\chi(\bar{M})<0$ the group $H_1(\bar{M})$ is not trivial, hence there are non-contractible edges of $\mathcal{N}(f)$.

Recall that the Bers theorem provides the following expresion for an eigenfunction $f$.
\begin{equation}
f(x)=r^n(A\sin (n\theta)+B\cos (n\theta))+O(r^{n+1}). \label{b}
\end{equation}
In the case $\dim (U_n)=2$ the ratio $(A{:}B)$ given by the expression~\eqref{b} is a homogeneous coordinate on both $P(U_n)$ and $F_x(2n)$. This implies that if $\dim (U_n)=2$ then the mapping $f \mapsto s(\mathcal{N}(f),x)$ is a diffeomorphism between $U_n$ and $F_x(2n)$.  In particular, it sends the generator of $\pi_1(P(U_n))=\pi_1(S^1)=\mathbb{Z}$ into a generator of $\pi_1(F_x(2n))$, which is not contractible. That contradicts with Lemma 2, so the case $\dim (U_n)=2$ is excluded. This gives us the smooth section $\sigma(x)$ due to Proposition 1. Now we extend it onto $\bar{M}$.

\begin{lema}
Let ${\bar{M}}$ be endowed with a smooth metric coinciding with the one on $M$ everywhere except $Op((\partial M)^i)$. Then the section $\sigma$ can be extended from $M\setminus Op(\bigcup (\partial M)^i)$ to a section $\bar{\sigma}$ of $E_{\bar{M}}$ over the whole $\bar{M}$.
\end{lema}

\begin{lema}
Let $N$ be a closed surface. If there is a section $s$ of $E_N(2n)$, then $\chi(N)=0$.
\end{lema}

Proofs are postponed till the section~\ref{l34}.

\medskip

Applying Lemma 4 to the section $\bar{\sigma}$ constructed in Lemma 3 we obtain a contradiction since $\chi(\bar{M})\neq 0$ was supposed. Hence the assumption of Lemma 1 saying that Theorem 1 does not hold was wrong. $\square$

\section{Proof of Lemma 1}

\label{l1}

According to the Bers theorem any eigenfunction $f$ has the following leading term in polar coordinates $(r,\theta)$:
\begin{equation}
f(r,\theta)=r^k(A\sin(k\theta)+B\cos(k\theta))+O(r^{k+1}). \label{bx}
\end{equation}
Number $k$ is called an order of zero at $r=0$. There is a filtration $$U=U_0\supset U_1\supset \dots \supset U_k\supset \dots$$ of the spece $U$ where $U_k$ denotes the vector space consisting of functions $f\in U$ which have zero at $r=0$ of order at least $k$, or, equivalently, $\deg |_{\mathcal{N}(f)}(x)\geq 2k$. Each adjoint factor $U_k/U_{k+1}$ is at most 2-dimensional since the numbers $(A,B)$ from~\eqref{bx} specify the element of $U_k/U_{k+1}$ unambiguously. The first factor $U_0/U_1$ is at most 1-dimensional, since $\sin(0\theta)=0$. Summing the dimensions up we obtain the estimate $\dim (U_k)\geq \dim (U)-2k+1$ for $k>0$. It implies the first statement of Lemma 1 saying $\dim (U_n)\geq 1$ since and it was supposed that $\dim (U)\geq 2n$.

For any $f \in U_{n}$ consider a nodal graph $\mathcal{N}(f)$ in $\bar{M}$. Contract all the vertices along the edges towards the vertex $x$. Add new edges with endpoints at $x$ making discs from non-simply-connected faces. We have obtained a new graph $\mathcal{N}'$. If some of propositions of Lemma 1 fail, namely, there were other vertices except $x$, or there were non-simply-connected faces of $\mathcal{N}$, or $\deg|_{\mathcal{N}(f)} (x)\geq 2n$, or $\dim (U_n)>2$, then $\deg_{\mathcal{N}'}(x)\geq 2n+2$ and hence there is at least $n+1$ edges. The number of faces of $\mathcal{N}$ is no more than $l$ by the definition of $l$ and is the same for $\mathcal{N}'$. Now Euler characteristic of $\bar{M}$ can be estimated as 
$$\chi (\bar{M})=\#\{\mbox{vertices of }{\mathcal{N}'}\}-\#\{\mbox{edges of }{\mathcal{N}'}\}+\#\{\mbox{faces of }{\mathcal{N}'}\}\leq$$ $$\leq 1-(n+1)+l=-(l-\chi(\bar{M})+1)+l=\chi (\bar{M})-1$$
which is a contradiction. 

Finally, if it turned out that $\dim (U_n)>2$, then $\dim (U_{n+1})>0$ and there is an eigenfunction $f\in U_{n+1}$ whose nodal graph has the vertex $x$ of degree at least $2n+2$, which we excluded above.

Hence all the ways Lemma 1 could fail lead to the contradiction and Lemma 1 is completely proven.

\section{Proof of Lemma 2}

\label{l2}

First of all we need a technical proposition.

\begin{prop}
Suppose $f(p)$ is a continuous family of functions depending on a parameter $p$, each $f(p)$ satisfy $\Delta f(p)=\lambda f(p)$ for some real number $\lambda \in \mathbb{R}$ and the nodal graph $\mathcal{N}(p)= \mathcal{N}(f(p))$ has only one vertex $x(p)$ in $\bar{M}$ and $\deg_{\mathcal{N}(p)}(x(p))=2n$. Then the loops provided by edges of $\mathcal{N}(p)$ do not change their homotopy class in the local system $\{\pi_1(\bar{M},x(p))\}$ while $p$ changes.
\end{prop}

{\bfseries Proof}. As in the proof of the Proposition 1, for any fixed value $p$ of the parameter consider the function $f(p)$ on $\bar{M}$ and take $n$-th term $f^{(n)}$ of its Taylor series at $x(p)$. Then according to Bers theorem we have the following formula
$$f^{(n)}(p)=r^n_x(A(p)\sin(n\theta_x)+B(p)\cos(n\theta_x))+O(r^{n+1}_x)$$
for some polar coordinates $(r_x,\theta_x)$ centered at $x(p)$. Note that $(A(p),B(p))\neq 0$ for any $p$ since $\deg_{\mathcal{N}(p)} (x(p))=2n$. This implies by the implicit function theorem that near $x(p)$ the graph $\mathcal{N}(p)$ undergoes an isotopy while $p$ changes. For pieces of the graph in~$\bar{M}\setminus Op(\bigcup (\partial M)^i \cup (x(p)))$ argumentation is analogous and even simplier.

It remains to find out what happens near $\bigcup pt^i$. We do not prove here that $\mathcal{N}$ undergoes an isotopy there as well. The invariance of homotopical class of loops of $\mathcal{N}(p)$ is provided by the following observation. If at some value $p_0$ one gets $pt^i\in \mathcal{N}(p_0)$ then the restriction of $f(p)$ on $\partial Op(pt^i)\cong S^1$ has two zeroes $\partial Op(pt^i)\cap \mathcal{N}(p_0)={:}\{in^i,out^i\}$ and they are non-critical. Then there is exactly two non-critical zeroes $\{in^i(p),out^i(p)\}$ of $f(p)|_{\partial Op(pt^i)}$ for $p\in Op(p_0)$ as well and $\{in^i(p),out^i(p)\}$ depends continuously on $p$. Recalling that $\mathcal{N}(f(p))$ has no vertices except $x$ and $x\notin Op(pt^i)$, we obtain that for $p \in Op(p_0)$ the interesting piece $\mathcal{N}(p)\cap Op(pt^i)$ has one connected component, no vertices and precisely two endpoints $\{in^i(p),out^i(p)\}$. That means $\mathcal{N}(p)\cap Op(pt^i)$ is an arc joining $\{in^i(p),out^i(p)\}$ inside $Op(pt^i)$ for $p \in Op(p_0)$. 

Finally, for any value $p_0$ of the parameter we can split any arc $l\subset \mathcal{N}(p)$ into pieces $l\cap Op(pt^i)$ (for those $i$ whose $pt^i$ are contained in $\mathcal{N}(p_0)$) and the rest \lq regular\rq\ part of $l$. We showed above that the \lq regular\rq\ part undergoes an isotopy for $p\in Op(p_0)$, the we showed that the pieces $l\cap Op(pt^i)$ preserve also their classes in local systems $\{\pi_1(\bar{M},in^i(p),out^i(p))\}$. That implies that the class of $l$ itself is preserved.$\square$

\medskip

Now let us prove Lemma 2. First of all, an isotopy of a graph does not change the class of its loops in $\{\pi_1(\bar{M},x(t))\}$. Hense Proposition 3 is valid not only for nodal graphs but for any graphs isotopic to nodal ones as well. Let us refer to the loop $s(\mathcal{N}(t),x)\in \pi_1(F_x)$ as a \textit{turn} of a star. If the turn of the star at a point $x$ is not null-homotopic, then its $2n$-multiple induces a turn by $2\pi m$ ($m\in \mathbb{Z}\setminus \{ 0\}$), which implies that every ray at $x$ comes back to its initial position. From now on consider not the initial isotopy, but this its multiple one.

Let us give the sketch of the next part of the proof of Lemma 2. When the star turns at $2\pi m$, every loop is conjugated by the loop $\gamma^m$ where $\gamma$ generates $\pi_1(Op(x)\setminus \{ x\})$. Although $\gamma$ is contractible in $\bar{M}$, we can use the information about the direction the ends of the loop go from. We will restrict to a class of loops ``remembering those directions''. In this class the loop $\gamma$ will be non-contractible and, moreover, its multiples will commute only with themselves. Recall that Lemma 2 assumes that there is a non-contractible edge $l$ of $\mathcal{N}$. Then the loops $\gamma^m$ and $l$ do not commute. But on another hand the $s^1$-isotopy assumed by Lemma 2 we have $\gamma^ml\gamma^{-m}=l$. This contradiction shows that the assumption $s(\mathcal{N}(t),x)\neq 0$ in $\pi_1(F_x)$ was wrong.

Now we produce a formal proof. 

Take polar coordinates $(r,\theta)$ centered at $x$ and take $\varepsilon >0$ such that in the ball $B_\varepsilon:\{(r,\theta)|r\leq\varepsilon\}$ the nodal lines differ slightly from straight rays $\theta=const$ (it is enough that $\frac{\partial r}{\partial \xi}>0$ for some parameter $\xi$ on a ray). Let our family $\mathcal{N}$ be parameterized by $t\in S^1$. Now we freeze it near $x$, i.e. consider another homotopy $\mathcal{N}'(t)$ which coincides with $\mathcal{N}(t)$ outside $B_\varepsilon$ and defined inside as $\mathcal{N}'(t)=\mathcal{N}(t\frac{2r-1}{\varepsilon})$ on $B_\varepsilon\setminus B_\frac{\varepsilon}{2}$ and $\mathcal{N}'(t)$ be constant inside $B_\frac{\varepsilon}{2}$. Graphs $\mathcal{N}'(0)$ and $\mathcal{N}'(1)$ differ only inside $B_\varepsilon$. The rays of the first one go almost straight along radii when the rays of the second one go $m$ times around $x$ before leaving $B_\varepsilon$.

Take now the non-contractible loop $l$ of $\mathcal{N}'(0)$. We can throw away its fixed end parts from inside $B_\frac{\varepsilon}{2}$ and join the cuts by an arc ($t$-independent) going along $r=\frac{\varepsilon}{2}$. We obtain a homotopy of a loop $l'(t)$ in $\bar{M}\setminus B_\frac{\varepsilon}{2}=M_B$, such that $l'(0)$ is non-contractible in $\bar{M}$ (as $l$ was) and $l'(1)=\gamma^ml'\gamma^{-m}$ where $\gamma$ generates $\pi_1(\{(\frac{\varepsilon}{2},\theta)\})$. So, $\gamma^m$ commutes with $l'$. To obtain a contradiction we need the following proposition.

\begin{prop}
Let $N$ be a closed surface $\chi (\bar{N})<0$, $B$ be a disc in $N$ and $\gamma$ be the image in $\pi_1 (N\setminus B)$ of a generator of $\pi_1 (\partial B)$. Then the centralisator $Z(\gamma)\subset \pi_1 (N\setminus B)$ is the group $\left\langle \gamma \right\rangle$ generated by $\gamma$.
\end{prop}

{\bfseries Proof}. Denote $N\setminus B$ by $N_B$. Take the cover $pr: \tilde{N}_B \rightarrow N_B$ induced by universal cover $\tilde{N}$ of $N$. The surface $\tilde{N}_B$ is obtained from a disc $\tilde{N}$ by cutting holes $pr^{-1}(B)$ numbered by elements $e_i$ of $\pi_1(N)$. Take non-intersecting paths $\gamma_{e_i}$ joining the marked point of $\tilde{N}_B$ with its image under the action of $e_i$. Then every element of $\pi_1(\tilde{N}_B$ can be represented uniquely as a product of terms $e_i\circ e_i(\gamma^{\pm 1})\circ e_i^{-1}$ with maybe one $e_j$ at the end where no two adjacent terms cancels. From such representation it is clear that $\gamma^m$ for $m\neq 0$ commute only with its multiple. $\square$

Applying Proposition 4 to the loop $l'$ we deduce that $l'$ is a power of $\gamma$ and hence is contractible in $\bar{M}$. But it was taken as a loop homotopic to non-contractible $l$. This contradiction shows the initial assumption of non-zero turn was false, which proves Lemma 2. $\square$

\section{Proofs of Lemmas 3 and 4}

\label{l34}
{~}

Proof of Lemma 3. 

We extend the section $\sigma$ of $E_M(2n)$ from $M\setminus \bigcup(D^i)$ to the section $\bar{\sigma}$ of $E_{\bar{M}}(2n)$ over the whole $\bar{M}$ where $D^i$ are disc neighborhoods of contracted holes $(\partial M)^i$. The corresponding obstruction $h$ to the extension $\sigma$ on $D^i$ lies in $$H^2(D^i,\partial(D^i);\pi_1(F_x))=\mathbb{Z}.$$
Here the Lemma 2 would be useful to conclude that the obstruction vanishes. But the obstruction is given by $s(\mathcal{N}(f_{x(t)}),x(t))$ with the vertex $x(t)$ running along $\partial(D^i)$, while Lemma 2 deals only with a fixed vertex $x$. So we only have to put all vertices of $\mathcal{N}(f_x(t))$ to the same point. 

Take a continuous family $R_{x,\tau}$ of diffeomorphisms of $Op(D^i)$ ($x\in Op(D^i),\tau\in [0,1]$), such that $R_{x,0}=R|_{\partial Op(D^i)}=\mbox{\upshape id}$ and $R_{x,1}(l(x))=pt^i=(\partial M)^i/(\partial M)^i$. Obstruction $h$ equals to the loop given by the section $s(\mathcal{N}(f_{x(t)}),x(t))$ in $\pi_1(F(2n))$ for the standard trivialisation of $E_{D^i}(2n)$. Now we can move it by $R_{\partial(Op(D^i)),1}$ to the single fibre $F_{pt^i}$. We can apply Lemma 2 to the graphs $\mathcal{N}'(t)$ obtained by isotopy $R_{x,\tau}$ from $\mathcal{N}_{f_x(t)}$ we can apply Lemma 2. Then the turn of their star $s(\mathcal{N}'(t),pt^i)$ is homotopically zero, so is the obstruction $h$. It means that we can extend the section $\sigma$ of $E_{\bar{M} \setminus \bigcup D^i}(2n)$ to the neighborhood $D^i$ of every $pt^i$ and hence to the whole $\bar{M}$. $\square$

\medskip
Proof of Lemma 4. 

First of all, lift the section $s$ from $N$ to its orienting cover $\tilde{N}$. For the oriented surface $\tilde{N}$ we have $E_{\tilde{N}}(2n)=S(T\tilde{N})/\rho_{\frac{\pi}{n}}$ where $\rho_{\frac{\pi}{n}}$ is the rotation in the spherisation $S(T\tilde{N})$ by ${\frac{\pi}{n}}$ in the positive direction. Then the Euler class $e(E_{\tilde{N}}(2n))$ equals $(2n)e(ST\tilde{N})$. The first fibration $E_{\tilde{N}}(2n)$ admits the section $s$, hence $e(E_{\tilde{N}}(2n))=0$. This implies $e(ST\tilde{N})=e(E_{\tilde{N}}(2n))/2n=0$ and $\chi(N)=\chi(\tilde{N})/2=0$. $\square$

\section*{Acknowledgements}
The author thanks A.V. Penskoi for many useful discussions and invaluable help in preparing this manuscript.

The article is supported in part by the Simons Foundation.

\end{document}